\documentclass{article}
\usepackage{amssymb}

\newtheorem{theorem}{Theorem}[section]

\newtheorem{corollary}[theorem]{Corollary}

\newtheorem{proposition}[theorem]{Proposition}

\newenvironment{proof}[1][Proof]{\textbf{#1.} }{\ \rule{0.5em}{0.5em}}

\begin{document}

\title{Geometric properties of Lagrangian mechanical systems}
\author{\textsc{Ioan Bucataru, Radu Miron} \\
{\small Faculty of Mathematics, ``Al.I.Cuza'' University,} \\
{\small Ia\c{s}i, 700506, Romania,} \\
{\small bucataru@uaic.ro, radu.miron@uaic.ro}}
\date{}
\maketitle

\begin{abstract}
The geometry of a Lagrangian mechanical system is determined by its
associated evolution semispray. We uniquely determine this semispray
using the symplectic structure and the energy of the Lagrange space
and the external force field. We study the variation of the energy
and Lagrangian functions along the evolution and the horizontal
curves and give conditions by which these variations vanish. We
provide examples of mechanical systems which are dissipative and for
which the evolution nonlinear connection is either metric or
symplectic.
\end{abstract}



\section*{Introduction}

A geometric approach of Riemannian mechanical systems has been
proposed recently by Munoz-Lecanda and Yaniz-Fernandez in
\cite{[munoz]}. Using techniques that are specific to Lagrange
geometry, R.Miron, \cite{[miron1]}, introduced and investigated some
geometric aspects of Finslerian and Lagrangian mechanical systems.
In this work we extend such geometric investigation of Lagrangian
mechanical systems. We determine the evolution semispray of a
mechanical system by using the symplectic structure and the energy
of the associated Lagrangian function and the external force field.

If the Lagrangian function is not homogeneous of second degree with
respect to the velocity-coordinates, as it happens in the Riemannian
and Finslerian framework, the energy of the system is different from
the Lagrangian function and the evolution curves (solution of the
Euler-Lagrange equations) are different from the horizontal curves
of the system. In this paper we study the variation of both energy
and Lagrangian function along the evolution curves and horizontal
curves. As it has been shown for the Riemannian case,
\cite{[munoz]}, we prove that the energy is decreasing along the
evolution curves of the system if and only if the external force
field is dissipative.

The canonical nonlinear connection of a Lagrange manifold is the
unique nonlinear connection that is metric and symplectic, as it has
been shown in \cite{[bucataru1]}. Conditions by which the evolution
nonlinear connection is either metric or symplectic are determined
in terms of the symmetric or skew-symmetric part of a (1,1)-type
tensor field associated with the external force field.

In the last part of the paper a special attention is paid to the
particular case of Finslerian mechanical systems. Examples of
dissipative mechanical systems are given.

\section{Geometric structures on tangent bundle}

In this section we introduce the geometric structures that live on
the total space of tangent (cotangent) bundle, which we are going to
use in this work such as: Liouville vector field, semispray,
vertical and horizontal distribution.

For an $n$-dimensional $C^{\infty}$-manifold $M$, we denote by $(TM,
\pi, M)$ its tangent bundle and by $(T^*M, \tau, M)$ its cotangent
bundle. The total space $TM$ ($T^*M$) of the tangent (cotangent)
bundle will be the phase space of the coordinate velocities
(momenta) of our mechanical system. Let $(U, \phi=(x^i))$ be a local
chart at some point $q\in M$ from a fixed atlas of
$C^{\infty}$-class of the differentiable manifold $M$. We denote by
$(\pi^{-1}(U), \Phi=(x^i, y^i))$ the induced local chart at $u\in
\pi^{-1}(q) \subset TM$. The linear map $\pi_{*,u}:T_uTM \rightarrow
T_{\pi(u)}M$ induced by the canonical submersion $\pi$ is an
epimorphism of linear spaces for each $u\in TM$. Therefore, its
kernel determines a regular, $n$-dimensional, integrable
distribution $V:u\in TM \mapsto V_uTM:=\mathrm{Ker} \pi_{*,u}
\subset T_uTM$, which is called the \textit{vertical distribution}.
For every $u \in TM$, $\{{\partial}/{\partial y^i}|_u\}$ is a basis
of $V_uTM$, where $\{{\partial}/{\partial x^i}|_u,
{\partial}/{\partial y^i}|_u\}$ is the natural basis of $T_uTM$
induced by a local chart. Denote by $\mathcal{F}(TM)$ the ring of
real-valued functions over $TM$ and by $\mathcal{X}(TM)$ the
$\mathcal{F}(TM)$-module of vector fields on $TM$. We also consider
$\mathcal{X}^v(TM)$ the $\mathcal{F}(TM)$-module of vertical vector
fields on $TM$. An important vertical vector field is
$\mathbb{C}=y^i({\partial}/{\partial y^i})$, which is called the
\textit{Liouville vector field}.

The mapping $J: \mathcal{X}(TM) \rightarrow \mathcal{X}(TM)$ given
by $J=({\partial}/{\partial y^i}) \otimes dx^i $ is called the
\textit{tangent structure} and it has the following properties: Ker
$J$ = Im $J$ = $\mathcal{X}^v(TM)$; rank $J=n$ and $J^2=0$. One can
consider also the cotangent structure $J^*=dx^i \otimes
({\partial}/{\partial y^i})$ with similar properties.

A vector field $S\in\chi(TM)$ is called a \textit{semispray}, or a
second order vector field, if $JS=\mathbb{C}$. In local coordinates
a semispray can be represented as follows:
\begin{equation}
S=y^i\frac{\partial}{\partial x^i}-2G^i(x,y)
\frac{\partial}{\partial y^i}. \label{smispray}
\end{equation}
Integral curves of a semispray $S$ are solutions of the following
system of SODE:
\begin{equation}
\frac{d^2x^i}{dt^2}+2G^i\left(x,\frac{dx}{dt}\right)=0. \label{sode}
\end{equation}
A \textit{nonlinear connection} $N$ on $TM$ is an $n$-dimensional
distribution $N: u\in TM \mapsto N_uTM\subset T_uTM$ that is
supplementary to the vertical distribution. This means that for
every $u\in TM$ we have the direct sum
\begin{equation}
T_uTM=N_uTM \oplus V_uTM. \label{thv} \end{equation} The
distribution induced by a nonlinear connection is called the
\textit{horizontal distribution}. We denote by $h$ and $v$ the
horizontal and the vertical projectors that correspond to the above
decomposition and by $\mathcal{X}^h(TM)$ the
$\mathcal{F}(TM)$-module of horizontal vector fields on $TM$. For
every $u=(x,y)\in TM$ we denote by ${\delta}/{\delta
x^i}|_u=h({\partial}/{\partial x^i}|_u).$ Then $\{{\delta}/{\delta
x^i}|_u, {\partial}/{\partial y^i}|_u\}$ is a basis of $T_uTM$
adapted to the decomposition (\ref{thv}). With respect to the
natural basis $\{{\partial}/{\partial x^i}|_u, {\partial}/{\partial
y^i}|_u\}$ of $T_uTM$, we have the expression:
\begin{equation}
\left.\frac{\delta}{\delta x^i}\right|_u =
\left.\frac{\partial}{\partial x^i}\right|_u -
N^j_i(u)\left.\frac{\partial}{\partial y^j}\right|_u, \ u \in TM.
\label{coefficients}\end{equation} The functions $N^i_j(x,y)$,
defined on domains of induced local charts, are called the
\textit{local coefficients} of the nonlinear connection. The
corresponding dual basis is $\{dx^i, \delta y^i=dy^i + N^i_jdx^j\}$.

It has been shown by M. Crampin \cite{[crampin1]} and J. Grifone
\cite{[grifone]} that every semispray determines a nonlinear
connection. The horizontal projector $h$ that corresponds to this
nonlinear connection is given by:
\begin{equation}
h(X)=\frac{1}{2}\left(Id-\mathcal{L}_SJ\right)(X)=\frac{1}{2}\left(
X-[S,JX] - J[S,X]\right). \label{hproj} \end{equation} Local
coefficients of the induced nonlinear connection are given by
$N^i_j={\partial G^i}/{\partial y^j}$.

\section{Geometric structures on a Lagrange space}

The presence of a regular Lagrangian on the tangent bundle $TM$
determines the existence of some geometric structures one can
associate to it such as: semispray, nonlinear connection and
symplectic structure.

Consider $L^{n}=(M, L)$ a Lagrange space. This means that $L:TM
\longrightarrow \mathbb{R}$ is differentiable of $C^{\infty}$-class
on $\widetilde{TM}=TM\setminus \{0\}$ and only continuous on the
null section. We also assume that $L$ is a regular Lagrangian. In
other words, the (0,2)-type, symmetric, d-tensor field with
components
\begin{equation} g_{ij}=\frac{1}{2}\frac{\partial^2 L}{\partial
y^i\partial y^j} \textrm{ has rank n on } \widetilde{TM}.
\label{lmetric}
\end{equation}
The Cartan 1-form $\theta_L$ of the Lagrange space can be defined as
follows:
\begin{equation}
\theta_L=J^*(d L)=d_JL=\frac{\partial L}{\partial y^i} dx^i.
\label{theta}
\end{equation}
For a vector field $X=X^i({\partial}/{\partial x^i}) +
Y^i({\partial}/{\partial y^i})$ on $TM$, the following formulae are
true:
\begin{equation}
\theta_L(X)=dL(JX)=d_JL(X)=(JX)(L)=\frac{\partial L}{\partial
y^i}X^i. \label{thetax}
\end{equation}
The Cartan 2-form $\omega_L$ of the Lagrange space can be defined as
follows:
\begin{equation}
\omega_L=d \theta_L = d(J^*(dL))=dd_JL= d\left(\frac{\partial
L}{\partial y^i}dx^i\right). \label{omega} \end{equation} In local
coordinates, the Cartan 2-form $\omega_L$ has the following
expression:
\begin{equation}
\omega_L =2g_{ij} dy^j \wedge d x^i +
\displaystyle\frac{1}{2}\left(\displaystyle\frac{\partial^2
L}{\partial y^i\partial x^j} - \displaystyle\frac{\partial^2
L}{\partial x^i \partial y^j}\right)d x^j\wedge d x^i.
\label{omegal}
\end{equation}
We can see from expression (\ref{omegal}) that the regularity of the
Lagrangian $L$ is equivalent with the fact that the Cartan 2-form
$\omega_L$ has rank 2n on $\widetilde{TM}$ and hence it is a
symplectic structure on $\widetilde{TM}$. With respect to this
symplectic structure, the vertical subbundle is a Lagrangian
subbundle of the tangent bundle.

The canonical semispray of the Lagrange space $L^n$ is the unique
vector field $\mathring{S}$ on $TM$ that satisfies the equation
\begin{equation}
i_{\mathring{S}}\omega_L = -dE_L. \label{isomega}
\end{equation}
Here $E_L=\mathbb{C}(L)-L$ is the energy of the lagrange space
$L^n$. The local coefficients $\mathring{G}{}^i$ of the canonical
semispray $\mathring{S}$ are given by the following formula:
\begin{equation}
\mathring{G}{}^i=\frac{1}{4}g^{ik}\left(\frac{\partial^2 L}{\partial
y^k
\partial x^h} y^h-\frac{\partial L}{\partial x^k}\right).
\label{gi}
\end{equation}
Using the canonical semispray $\mathring{S}$ we can associate to a
regular Lagrangian $L$ a canonical nonlinear connection with local
coefficients given by expression $\mathring{N}{}_j^i={\partial
\mathring{G}{}^i}/{\partial y^j}$.

The horizontal subbundle $NTM$ that corresponds to the canonical
nonlinear connection is a Lagrangian subbundle of the tangent bundle
$TTM$ with respect to the symplectic structure $\omega_L$. This
means that $\omega_L(hX,hY)=0$, $\forall X,Y \in \chi(TM)$. For a
more detailed discussion regarding symplectic structures in Lagrange
geometry we recommend \cite{[anastasiei]}. In local coordinates this
implies the following expression for the symplectic structure
$\omega_L$:
\begin{equation}
\omega_L=2g_{ij}\mathring{\delta}y^j \wedge dx^i. \label{omegalocal}
\end{equation}

The dynamical derivative that corresponds to the pair
$(\mathring{S}, \mathring{N})$ is defined by $\nabla: \chi^v(TM)
\longrightarrow \chi^v(TM)$ through:
\begin{equation}
\nabla \left(X^i\frac{\partial}{\partial
y^i}\right)=\left(\mathring{S}(X^i)+ X^j\mathring{N}{}^i_j\right)
\frac{\partial}{\partial y^i}. \label{covariant}
\end{equation}
In terms of the natural basis of the vertical distribution we have
\begin{equation}
\nabla \left(\frac{\partial}{\partial y^i}\right)=\mathring{N}{}^j_i
\frac{\partial}{\partial y^j}. \label{ncovariant}
\end{equation}
Hence, $\mathring{N}{}^i_j$ are also local coefficients of the
dynamical derivative. Dynamical derivative $\nabla$ is the same with
the covariant derivative $D$ in \cite{[crampin2]} or
$\mathcal{D}_{\Gamma}$ in \cite{[krupkova1]}, where it is called the
$\Gamma$-derivative. Dynamical derivative $\nabla$ has the following
properties:
\begin{itemize}
\item[1)] $\nabla(X+Y)=\nabla X + \nabla Y, \forall X, Y \in
\chi^v(TM)$, \vspace*{1mm} \item[2)] $\nabla(fX)=S(f) X + f\nabla X,
\forall X \in \chi^v(TM), \forall f\in \mathcal{F}(TM)$.
\end{itemize}
It is easy to extend the action of $\nabla$ to the algebra of
d-tensor fields by requiring for $\nabla$ to preserve the tensor
product. For the metric tensor $g$, its dynamical derivative is
given by
\begin{equation}
(\nabla g)(X,Y)=S(g(X,Y)) - g(\nabla X,Y) - g(X,\nabla Y), \forall
X,Y. \label{covariantg}
\end{equation}
In local coordinates, we have:
\begin{equation}
g_{ij|}:=(\nabla g)\left(\frac{\partial}{\partial y^i},
\frac{\partial}{\partial y^j}\right)=\mathring{S}(g_{ij}) -
g_{im}\mathring{N}{}^m_j - g_{mj}\mathring{N}{}^m_i.
\label{covariantgij}
\end{equation}
The canonical nonlinear connection $\mathring{N}$ is metric, which
means that $\nabla g=0$. In \cite{[bucataru1]} it is shown that the
canonical nonlinear connection of a Lagrange space is the unique
nonlinear connection that is metric and symplectic.

\section{Geometric structures of Lagrangian mechanical systems}

The dynamical system of a Lagrangian mechanical system is a
semispray, which we call the evolution semispray of the system. Such
semispray is uniquely determined by the symplectic structure and the
energy of the underlying Lagrange space and the external force
field. The energy of the system is decreasing if and only if the
force field is dissipative.

The nonlinear connection we associate to the evolution semispray is
called the evolution nonlinear connection. Conditions by which such
nonlinear connection is either metric or symplectic are studied.

A \textit{Lagrangian mechanical system} is a triple
$\Sigma_L=(M,L,V)$, where $(M,L)$ is a Lagrange space and
$V=V^i(x,y)({\partial}/{\partial y^i})$ is a vertical vector field,
which is called the external force field of the system. Using the
metric tensor $g_{ij}$ of the Lagrange space one can define the
vertical one-form $\sigma=\sigma_i dx^i$, where
$\sigma_i(x,y)=g_{ij}(x,y)V^j(x,y)$.

The external force field $V$ is \textit{dissipative} if $g({\mathbb
C},V)=g_{ij}y^i V^j\leq 0$.

The evolution equations of the mechanical system $\Sigma_L$ are
given by the following Lagrange equations, \cite{[miron1]}:
\begin{equation}
\frac{d}{dt}\left(\frac{\partial L}{\partial y^i}\right) -
\frac{\partial L}{\partial x^i}=\sigma_i, \ y^i=\frac{dx^i}{dt}.
\label{eqev1} \end{equation} For a regular Lagrangian, Lagrange
equations (\ref{eqev1}) are equivalent to the following system of
second order differential equations:
\begin{equation}
\frac{d^2x^i}{dt^2} + 2\mathring{G}{}^i\left(x,
\frac{dx}{dt}\right)=\frac{1}{2}V^i\left(x,\frac{dx}{dt}\right).
\label{eqev2}\end{equation} Here the functions $\mathring{G}{}^i$
are local coefficients of the canonical semispray of the Lagrange
space, given by expression (\ref{gi}). Since $(1/2)V^i$ are
components of a d-vector field on $TM$, from expression
(\ref{eqev2}) we obtain that the functions
\begin{equation}
2G^i(x,y)=2\mathring{G}{}^i(x,y)-\frac{1}{2}V^i(x,y) \label{evsem1}
\end{equation} are coefficients of a semispray $S$, to which we
refer to as the evolution semispray of the system. Therefore, the
evolution semispray is given by $S=\mathring{S} +(1/2)V$. Integral
curves of the evolution semispray $S$ are solutions of the SODE
given by expression (\ref{eqev2}).
\begin{theorem}
The evolution semispray $S$ is the unique vector field on $TM$,
solution of the equation
\begin{equation}
i_S\omega_L=-dE_L+\sigma. \label{evsem2} \end{equation}
\end{theorem}
\begin{proof}
Since $\omega_L$ is a symplectic structure on $TM$, equation
(\ref{evsem2}) uniquely determine a vector field $S$ on $TM$. The
vector field $V^i({\partial}/{\partial y^i})$ is the unique vector
field that satisfies $i_{V^i({\partial}/{\partial
y^i})}\omega_L=2\sigma$. Using the linearity of equations
(\ref{isomega}) and (\ref{evsem2}) we can see that $S=\mathring{S}+
(1/2) V^i({\partial}/{\partial y^i})$ is the unique solution of
equation(\ref{evsem2}).
\end{proof}
\begin{corollary} The energy of the Lagrange space $L^n$ is
decreasing along the evolution curves of the mechanical system if
and only if the external force field is dissipative.
\end{corollary}
\begin{proof} Using the fact that the evolution semispray $S$ is solution of
equation (\ref{evsem2}) and the skew-symmetry of $\omega_L$ we
obtain $S(E_L)=dE_L(S)=\sigma(S)=\sigma_iy^i$. Along the evolution
curves of the mechanical system, one can write this expression as
follows:
\begin{equation}
\frac{d}{dt}\left(E_L\right)=\sigma_i\left(x,\frac{dx}{dt}\right)\frac{dx^i}{dt}.
\label{dissip1} \end{equation} Therefore, the energy is decreasing
along the evolution curves if and only if $\sigma_iy^i\leq 0$.
\end{proof}

Expression (\ref{dissip1}) has been obtained in \cite{[munoz]} for
the particular case of a Riemannian mechanical system. For the
general case of Lagrangian mechanical system it has been obtained
also in \cite{[miron1]}, using different techniques.

The evolution nonlinear connection of the mechanical system
$\Sigma_L$ has the local coefficients $N^i_j$ given by
\begin{equation}
N^i_j=\frac{\partial G^i}{\partial y^j}=\frac{\partial
\mathring{G}{}^i}{\partial y^j} - \frac{1}{4}\frac{\partial
V^i}{\partial y^j} = \mathring{N}{}^i_j - \frac{1}{4}\frac{\partial
V^i}{\partial y^j}. \label{evconnection}
\end{equation}

\begin{theorem} \label{thm:metricconev}
The evolution nonlinear connection $N$ is metric if and only if the
$(0,2)$-type d-tensor field ${\partial \sigma_i}/{\partial y^j}$ is
skew-symmetric.\end{theorem}
\begin{proof}
The evolution nonlinear connection is metric if and only if the
dynamical covariant derivative of the metric tensor $g_{ij}$ with
respect to the pair $(S,N)$ vanishes. This covariant derivative is
given by
\begin{equation}
\begin{array}{ll}
g_{ij|} & =S(g_{ij})-g_{ik}N^k_j- g_{kj}N^k_i \vspace{2mm} \\
& = \left(\mathring{S} + \displaystyle\frac{1}{2}V \right)( g_{ij})
- g_{ik}\left(\mathring{N}{}^k_j -
\displaystyle\frac{1}{4}\frac{\partial V^k}{\partial y^j} \right) -
g_{kj}\left(\mathring{N}{}^k_i -
\displaystyle\frac{1}{4}\frac{\partial V^k}{\partial y^i}
\right) \vspace{2mm} \\
& = \mathring{S}(g_{ij})-g_{ik}\mathring{N}{}^k_j-
g_{kj}\mathring{N}{}^k_i + \displaystyle\frac{1}{4}\left( 2V(g_{ij})
+ g_{ik} \frac{\partial V^k}{\partial y^j} + g_{kj} \frac{\partial
V^k}{\partial y^i} \right)
\vspace{2mm}\\
& = \displaystyle\frac{V^k}{4}\left( 2\frac{\partial
g_{ij}}{\partial y^k} - \frac{\partial g_{ik}}{\partial y^j} -
\frac{\partial g_{kj}}{\partial y^i} \right) +
\displaystyle\frac{1}{4}\left(\frac{\partial \sigma_i}{\partial y^j}
+ \frac{\partial \sigma_j}{\partial y^i}\right)
\vspace{2mm} \\
& = \displaystyle\frac{1}{4}\left(\frac{\partial \sigma_i}{\partial
y^j} + \frac{\partial \sigma_j}{\partial y^i}\right).
\end{array} \label{eqmetricconev}
\end{equation}
In the above calculations we did use the fact that the canonical
nonlinear connection $\mathring{N}$ is metric and the Cartan tensor
\begin{equation}
C_{ijk}=\frac{1}{2}\frac{\partial g_{ij}}{\partial
y^k}=\frac{1}{4}\frac{\partial^3L}{\partial y^i\partial y^j\partial
y^k} \label{cartan} \end{equation} is totally symmetric.

The dynamical covariant derivative of the metric tensor $g_{ij}$
with respect to the pair $(S,N)$ is given by
\begin{equation}
g_{ij|}=\frac{1}{4}\left(\frac{\partial \sigma_i}{\partial y^j} +
\frac{\partial \sigma_j}{\partial y^i}\right).\label{dynconev}
\end{equation}
Consequently, the evolution nonlinear connection is metric if and
only if the $(0,2)$-type d-tensor field ${\partial
\sigma_i}/{\partial y^j}$ is skew-symmetric.
\end{proof}

\begin{theorem} \label{thm:symplecticconnev} The evolution nonlinear connection
is compatible with the symplectic structure if and only if the
$(0,2)$-type d-tensor field ${\partial \sigma_i}/{\partial y^j}$ is
symmetric.\end{theorem}
\begin{proof}
The evolution nonlinear connection $N$ is compatible with the
symplectic structure of the Lagrange space if and only if
$\omega_L(hX,hY)=0$, $\forall X, Y\in \chi(TM)$, where $h$ is the
corresponding horizontal projector.

Let us consider the almost symplectic structure:
\begin{equation}
\omega=2g_{ij}\delta y^j\wedge dx^i, \label{almostomega}
\end{equation}
with respect to which both horizontal and vertical subbundles are
Lagrangian subbundles. Using expressions (\ref{omegalocal}) and
(\ref{evconnection}) the canonical symplectic structure can be
expressed as follows:
\begin{equation}
\begin{array}{ll}
\omega_L & = 2g_{ij}\mathring{\delta}y^j\wedge dx^i =
2g_{ij}\left(\delta y^j + \displaystyle\frac{1}{4}\frac{\partial
V^j}{\partial y^k} dx^k\right) \wedge dx^i
\vspace{2mm}\\
& = 2g_{ij}\delta y^j\wedge dx^i  + \displaystyle\frac{1}{4} \left(
\frac{\partial \sigma_i}{\partial y^j} - \frac{\partial
\sigma_j}{\partial y^i}\right) dx^j \wedge dx^i \vspace{2mm}\\
& = \omega + \displaystyle\frac{1}{2}F_{ij}dx^j\wedge dx^i.
\end{array} \label{omegaev2}
\end{equation}
Here $F_{ij}$ is the helicoidal tensor of the mechanical system
$\Sigma_L$ introduced by R.Miron in \cite{[miron1]}:
\begin{equation}
F_{ij}=\frac{1}{2}\left(\frac{\partial \sigma_i}{\partial y^j} -
\frac{\partial \sigma_j}{\partial y^i}\right). \label{helicoidal}
\end{equation}

Therefore, the evolution nonlinear connection is compatible with the
symplectic structure if and only if $\omega_L=\omega$ which is
equivalent to the fact that the helicoidal tensor of the mechanical
system $\Sigma_L$, given by expression (\ref{helicoidal}), vanishes.
\end{proof}

\section{Variation of energy and Lagrangian functions}

In this section we study the variation of energy and Lagrangian
functions along the horizontal curves of the evolution nonlinear
connection.

For a Lagrangian mechanical system $\Sigma_L$ consider $S$ the
evolution semispray given by expression (\ref{evsem1}) and the
evolution nonlinear connection given by expression
(\ref{evconnection}). Let $h$ be the corresponding horizontal
distribution given by expression (\ref{hproj}). Then $hS$ is also a
semispray, its integral curves are called horizontal curves of the
evolution nonlinear connection. They are solutions of the following
system of SODE:
\begin{equation}
\nabla\left(\frac{dx^i}{dt}\right)= \frac{d^2x^i}{dt^2} +
N^i_j\left(x, \frac{dx}{dt}\right)\frac{dx^j}{dt} = 0.
\label{hcurves}
\end{equation}
In the previous section we studied the variation of the energy
function along the evolution curves of the system. We shall study
now, the variation of the energy and Lagrangian functions along the
horizontal curves (\ref{hcurves}). This way we can determine
external force fields such that the Lagrangian or the energy
functions are first integrals for the system (\ref{hcurves}).

\begin{theorem} \label{hvarL}
Consider $h$ the horizontal projector of the evolution nonlinear
connection. We have the following formula for the horizontal
differential operator $d_h$ of the Lagrangian $L$:
\begin{equation}
2d_hL=d_J(S(L))-\sigma. \label{dhl}
\end{equation}
In local coordinates, formula (\ref{dhl}) is equivalent with the
following expression for the horizontal covariant derivative of the
Lagrangian $L$:
\begin{equation}
2L_{|i}:=2\frac{\delta L}{\delta x^i} = \frac{\partial}{\partial
y^i}(S(L))-\sigma_i. \label{dhllocal}
\end{equation}
\end{theorem}
\begin{proof}
We prove first the following formulae regarding the Cartan 1-form
$\theta_L$ of the Lagrange space:
\begin{equation}
\mathcal{L}_S\theta_L = dL +\sigma. \label{ftheta}
\end{equation}
By differentiating $\iota_S\theta_L=\mathbb{C}(L)$ we obtain
$d\iota_S\theta_L=d\mathbb{C}(L)$. Using the expression of the Lie
derivative $\mathcal{L}_S=d\iota_S+\iota_Sd$ we obtain
$\mathcal{L}_S\theta_L=\iota_S\theta_L+d\mathbb{C}(L)$. From the
defining formulae (\ref{omega}) and (\ref{evsem2}) for $\omega_L$
and $S$ we obtain $\mathcal{L}_S\theta_L = dE_L+ \sigma+
d\mathbb{C}(L)=dL+\sigma$

In order to prove (\ref{dhl}) we have to show that for every $X\in
\chi(TM)$, we have that
$$ 2(d_hL)(X):=2dL(hX)= (JX)(S(L))-\sigma(X).$$
Using formula (\ref{ftheta}) we obtain
$$
\begin{array}{ll}
0 & =(\mathcal{L}_S\theta_L - dL- \sigma)(X) =
S\theta_L(X)-\theta_L[S,X]-dL(X) -\sigma(X) \vspace{2mm}\\
& = S((JX)(L))-J[S,X](L)-dL(X) -\sigma(X) \vspace{2mm}\\
& = [S,JX](L) + (JX)(S(L)) - J[S,X](L) - dL(X) -\sigma(X)\vspace{2mm}\\
& = (JX)(S(L)) - dL\left(X - [S,JX] + J[S,X]\right) -\sigma(X)\vspace{2mm}\\
& = (JX)(S(L))-\sigma(X)-dL(2hX) .
\end{array}
$$
Consequently, formula (\ref{dhl}) is true.

Due to the linearity of the operators involved in formula
(\ref{dhl}) we have that formulae (\ref{dhl}) and (\ref{dhllocal})
are equivalent.
\end{proof}
\begin{corollary} \label{corollary1}
If the external force field of the mechanical system satisfies the
equation:
\begin{equation}
\frac{\partial V^k}{\partial y^i}y^i\frac{\partial L}{\partial y^k}=
-2\mathbb{C}\left(\mathring{S}(L)\right) \label{hsl}
\end{equation}
then the Lagrangian $L$ is constant along the horizontal curves of
the evolution nonlinear connection.
\end{corollary}
\begin{proof} The Lagrangian $L$ is constant along the horizontal curves of the
evolution nonlinear connection if and only if $hS(L)=0$. If we
contract expression (\ref{dhllocal}) by $y^i$ we obtain
\begin{equation}
\begin{array}{ll}
2(hS)(L) = &2L_{|i}y^i=2\displaystyle\frac{\delta L}{\delta x^i}y^i
= \mathbb{C}(S(L))-\sigma(S) \vspace{2mm}\\
& =\mathbb{C}\left(\mathring{S}(L)\right)+ \displaystyle\frac{1}{2}
\displaystyle\frac{\partial V^k}{\partial y^i}
y^i\displaystyle\frac{\partial L}{\partial y^k}. \end{array}
\label{hslzero}
\end{equation}
Therefore we can see that $hS(L)=0$ if and only if the external
force field satisfies equation (\ref{hsl}).
\end{proof}

\begin{theorem} \label{hvarEL}
Consider $h$ the horizontal projector of the evolution nonlinear
connection. We have the following expression for the horizontal
differential operator $d_h$ of the energy $E_L$:
\begin{equation}
d_hE_L(X)=-\omega_L(\mathring{S}, hX). \label{dhel}
\end{equation}
In local coordinates, this is equivalent with the following
expression for the horizontal covariant derivative of the energy
$E_L$:
\begin{equation}
E_{L|i}:=\frac{\delta E_L}{\delta x^i} =
2g_{ij}\left(2\mathring{G}{}^j - \mathring{N}{}^j_ky^k\right) +
\frac{1}{2}g_{jk}\frac{\partial V^j}{\partial y^i}y^k.
\label{dhllocal1}
\end{equation}
\end{theorem}
\begin{proof}
We have that $d_hE_L(X)=(dE_L)(hX)=-\omega_L(\mathring{S},hX)$.
Since $h$ is the horizontal projector for the evolution nonlinear
connection, for a vector field $X=X^i({\partial}/{\partial x^i}) +
Y^i({\partial}/{\partial y^i})$ on $TM$ we have
$$ hX=X^i\frac{\delta}{\delta x^i}= X^i\frac{\mathring{\delta}}{\delta
x^i} + X^i\frac{1}{4}\frac{\partial V^j}{\partial
y^i}\frac{\partial}{\partial y^j}. $$ Using expression
(\ref{omegalocal}) for the symplectic structure $\omega_L$ we have
\begin{equation}
\begin{array}{ll}
d_hE_L(X) & =\left(-2g_{ij}\mathring{\delta} y^i\wedge
dx^j\right)(\mathring{S}, hX) \vspace{2mm} \\
& = 2g_{ij}X^i\left(2\mathring{G}{}^j - \mathring{N}{}^j_ky^k\right)
+ \displaystyle\frac{1}{2}g_{jk}\frac{\partial V^j}{\partial
y^i}y^kX^i, \end{array}\label{dhel1} \end{equation} and therefore we
proved both formulae (\ref{dhel}) and (\ref{dhllocal1}).
\end{proof}

\begin{corollary} \label{corollary2}
If the external force field of the mechanical system satisfies the
equation:
\begin{equation}
g_{jk}\frac{\partial V^j}{\partial y^i}y^ky^i = -4g_{ij}
\left(2\mathring{G}{}^j -\mathring{N}{}^j_ky^k\right)y^i
\label{hsel}
\end{equation}
then the energy $E_L$ is constant along the horizontal curves of the
evolution nonlinear connection.
\end{corollary}

\section{Finslerian mechanical systems} Finsler geometry corresponds
to the case when the Lagrangian function is second order homogeneous
with respect to the velocity coordinates. This has various
implications for the geometry of a Finsler space: the energy
coincides with the fundamental function of the space and it is
constant along the geodesic curves which are also horizontal curves
for the canonical nonlinear connection. Therefore the geometry of a
Finslerian mechanical system has some special features.

A Lagrange space $L^n=(M,L)$ reduces to a Finsler space $F^n=(M,F)$
if the Lagrangian function is second order homogeneous with respect
to the velocity coordinates. In this case we shall use the notation
$F^2(x,y)=L(x,y)$ and therefore by using Euler's theorem for
homogeneous functions we have:
$$ \mathbb{C}(F^2)=\frac{\partial F^2}{\partial y^i}y^i = 2F^2. $$
A first consequence of the homogeneity condition is that the energy
of a Finsler space coincides with the square of the fundamental
function of the space: $E_{F^2}=\mathbb{C}(F^2)-F^2=F^2$.

A Finslerian mechanical system is a triple $\Sigma_F=(M, F, V)$,
where $(M, F)$ is a Finsler space and
$V=V^i(x,y)({\partial}/{\partial y^i})$ is a vertical vector field.

The evolution equations of the Finslerian mechanical system are
given by Lagrange equations (\ref{eqev1}) where $L(x,y)=F^2(x,y)$,
which are equivalent with the system of second order differential
equations (\ref{eqev2}). The local coefficients $\mathring{G}{}^i$
of the canonical semispray $\mathring{S}$ of the Finsler space can
be written in this case as follows:
\begin{equation}
2\mathring{G}{}^i(x,y)=\gamma^i_{jk}(x,y)y^jy^k, \label{fspray}
\end{equation} where $\gamma^i_{jk}(x,y)$ are Chrystoffel symbols of
the metric tensor $g_{ij}(x,y)$. Therefore, the evolution semispray
$S$ of the mechanical system has the local coefficients $2G^i$ given
by
\begin{equation}
2G^i(x,y)=\gamma^i_{jk}(x,y)y^jy^k - \frac{1}{2}V^i(x,y).
\label{fmspray}
\end{equation}

The evolution curves of the mechanical system are solutions of the
following SODE:
\begin{equation}
\frac{d^2x^i}{dt^2} + \gamma^i_{jk}\left(x,\frac{dx}{dt}\right)
\frac{dx^j}{dt}\frac{dx^k}{dt}
-\frac{1}{2}V^i\left(x,\frac{dx}{dt}\right)=0. \label{fmev}
\end{equation}

Local coefficients $N^i_j$ of the evolution nonlinear connection are
given by expression (\ref{evconnection}). Due to the homogeneity
conditions one can express them as follows:
\begin{equation}
N^i_j(x,y)=\gamma^i_{kj}(x,y)y^k-\frac{1}{4}\frac{\partial
V^i}{\partial y^j}(x,y). \label{fmevcon} \end{equation} Therefore,
the horizontal curves of the evolution nonlinear connection are
solutions of the following SODE:
\begin{equation}
\frac{d^2x^i}{dt^2} + \gamma^i_{jk}\left(x,\frac{dx}{dt}\right)
\frac{dx^j}{dt}\frac{dx^k}{dt} -\frac{1}{4}\frac{\partial
V^i}{\partial y^j}\left(x,\frac{dx}{dt}\right)\frac{dx^j}{dt}=0.
\label{fmhcurv}
\end{equation}

\begin{proposition} Consider a Finslerian mechanical system $\Sigma_F=(M, F^2, V)$,
with the external force field $V$ homogeneous of order zero. Then
the energy function $F^2(x,y)$ is constant along the horizontal
curves of the evolution nonlinear connection.
\end{proposition}
\begin{proof}
If the external force field $V$ is homogeneous of order zero, then
by Euler theorem we have $({\partial V^i}/{\partial y^j})y^j=0$ and
the horizontal curves of the evolution nonlinear connection given by
expression (\ref{fmhcurv}) coincide with the geodesics of the
Finsler space. The energy function $F^2(x,y)$ is constant along the
geodesic curves of the Finsler space.

One can obtain this result by using also expressions (\ref{hsel}) or
(\ref{hsl}). The right hand side of both expression vanishes for a
Finsler space, while the left hand side vanishes due to the
homogeneity of the external force field.
\end{proof}

For a Finsler space $(M, F^2(x,y))$, the local coefficients
$2\mathring{G}{}^i(x,y)$ of the canonical semispray are given by
expression (\ref{fspray}) and therefore they are second order
homogeneous with respect to the velocity variables. This implies
that $2\mathring{G}{}^j = \mathring{N}{}^j_ky^k$ and equation
(\ref{dhllocal1}) can be written as follows:
\begin{equation}
F^2_{|i}=\frac{1}{2}g_{jk}\frac{\partial V^j}{\partial y^i}y^k=
\frac{1}{2}\frac{\partial \sigma_k}{\partial y^i}{y^k}. \label{f2i}
\end{equation}
Here we did use the symmetry of the Cartan tensor (\ref{cartan}) and
the zeroth homogeneity of the metric tensor $g_{ij}$. If the
helicoidal tensor of the mechanical system vanishes and the external
force field is zero homogeneous then the horizontal covariant
derivative of the energy function vanishes, in other words
$F^2_{|i}=0$.

\section{Examples} In this section we give
examples of Finslerian and Lagrangian mechanical systems which have
the properties that we studied in the previous sections.

\textbf{1.} Consider a Lagrangian mechanical system $\Sigma_L=(M, L,
V)$, where $V=e\mathbb{C}=ey^i({\partial}/{\partial y^i})$ and $e$
is a constant. We call this system a Liouville mechanical system.
The evolution semispray $S$ and nonlinear connection $N$ have the
local coefficients given by:
\begin{equation}
2G^i(x,y)=2\mathring{G}{}^i-\frac{e}{2}y^i, \
N^i_j=\mathring{N}{}^i_j -\frac{e}{4}\delta^i_j. \label{ex1}
\end{equation}

Therefore, the Liouville mechanical system $\Sigma_L=(M, L,
e\mathbb{C})$ has some special properties.

\begin{proposition} The helicoidal tensor $F_{ij}$ of the system vanishes.
Consequently, the evolution nonlinear connection is compatible with
the symplectic structure of the Lagrange space. \end{proposition}
\begin{proof}
The helicoidal tensor of the system is given by expression
(\ref{helicoidal}). Since $\sigma_i=eg_{ik}y^k$ we have the
following expression: \begin{equation} \frac{\partial
\sigma_i}{\partial y^j}= e\left( \frac{\partial g_{ik}}{\partial
y^j} y^k + g_{ij}\right) = e \left( 2C_{ijk}y^k + g_{ij}\right).
\label{ex:sigmaij} \end{equation} According to the above formula we
have that the $(0,2)$-type d-tensor field ${\partial
\sigma_i}/{\partial y^j}$ is symmetric and using Theorem
\ref{thm:symplecticconnev}, the evolution nonlinear connection is
compatible with the symplectic structure of the Lagrange space.
\end{proof}

\begin{proposition} The dynamical covariant derivative of the metric tensor $g_{ij}$
with respect to the pair $(S,N)$ is given by the following formula:
\begin{equation}
g_{ij|}=\frac{e}{2}\left(2C_{ijk}y^k+ g_{ij}\right). \label{ex:gij}
\end{equation} Consequently, the evolution nonlinear connection is metric
if and only if the metric tensor is homogeneous of order $-1$.
\end{proposition}
\begin{proof}
Formula (\ref{ex:gij}) follows immediately from expressions
(\ref{ex:sigmaij}) and (\ref{dynconev}). We have that $g_{ij|}=0$ if
and only if
$$ \frac{\partial g_{ij}}{\partial y^k}y^k =-1 \cdot g_{ij},$$ which is
equivalent with the homogeneity of order $-1$ of the metric tensor
$g_{ij}$.
\end{proof}

If the Liouville mechanical system is also a Finslerian one, then
$C_{ijk}y^k=0$ and consequently $g_{ij|}=(e/2)g_{ij}$.

\begin{proposition} We assume the Liouville mechanical system is Finslerian.
Then, the system is dissipative if and only if $e <0$.
\end{proposition}
\begin{proof}
The system is dissipative if and only if $0>g(\mathbb{C},
V)=eg_{ij}y^iy^j=eF^2$, which holds true if and only if $e<0$.
\end{proof}

The above result holds true even if the Liouville mechanical system
is truly Lagrangian, but in this case we have to ask for the
supplementary condition $g_{ij}(x,y)y^iy^j>0$. A sufficient
condition for this is the positive definiteness of the metric tensor
$g_{ij}$.

\textbf{2.} Consider the Finslerian mechanical system: $\Sigma_F=(M,
F^2(x,y), (e/F)\mathbb{C})$. For this system, the external force
field
$$ V=\frac{ey^i}{F}\frac{\partial}{\partial y^i} $$ is zero
homogeneous. According to the previous section the horizontal curves
of the evolution nonlinear connection coincide with the geodesic
curves of the Finsler space $(M, F^2)$.

\begin{proposition}
For the Finslerian mechanical system $\Sigma_F$ the following
properties hold true.
\begin{itemize}
\item[i)] The helicoidal tensor $F_{ij}$ of the system vanishes and
hence the evolution nonlinear connection is compatible with the
symplectic structure of the Finsler space.
\item[ii)] The horizontal covariant derivative of the energy
function vanishes. In other words $F^2_{|i}=0$ and hence $F^2$ is
constant along the horizontal curves of the evolution nonlinear
connection. \item[iii)] The system is dissipative if and only if
$e<0$. \end{itemize} \end{proposition}
\begin{proof}

i) Consider $y_i=g_{ij}y^j$. Using the symmetry of the Cartan tensor
and the zero homogeneity of the metric tensor $g_{ij}$ we obtain
${\partial y_i}/{\partial y^j}=g_{ij}$. Therefore we have:
\begin{equation}
\frac{\partial \sigma_i}{\partial y^j}= \frac{\partial}{\partial
y^j}\left(\frac{e y_i}{F}\right)= \frac{e^2
g_{ij}-\sigma_i\sigma_j}{e F}. \label{siyj} \end{equation} From
expression (\ref{siyj}) we obtain that the helicoidal tensor
$F_{ij}$ vanishes.

ii) The horizontal covariant derivative of $F^2$ is given by
expression (\ref{f2i}). Using the symmetry of the tensor ${\partial
\sigma_k}/{\partial y^i}$ and the zero homogeneity of the external
force field we obtain
\begin{equation}
F^2_{|i}=\frac{1}{2}\frac{\partial \sigma_k}{\partial y^i}{y^k} =
\frac{1}{2}\frac{\partial \sigma_i}{\partial y^k}{y^k}=0.
\label{f3i}
\end{equation}

iii) The system is dissipative if and only if $(e/F)g(\mathbb{C},
\mathbb{C})<0$, which is equivalent to $e<0$. \end{proof}

\vspace*{3mm}

An alternative approach of such systems can be obtained by means of
Legendre transformations by using the theory of Cartan and Hamilton
spaces. In this new framework, the Hamilton equations are used
instead of Euler-Lagrange equations and the Hamiltonian vector field
is used instead of the canonical semispray.

\vspace*{5mm}

\noindent\textbf{Acknowledgment} This work has been supported by
CEEX Grant 31742/2006 from the Romanian Ministry of Education.

\end{document}